\documentclass{amsart}
\usepackage{graphicx,amssymb}
\usepackage{a4wide}

\theoremstyle{plain}
\newtheorem{thm}{Theorem}[section]
\newtheorem{lem}[thm]{Lemma}
\newtheorem{prop}[thm]{Proposition}
\newtheorem{cor}[thm]{Corollary}

\theoremstyle{definition}
\newtheorem{defn}[thm]{Definition}
\newtheorem{exmp}[thm]{Example}

\def\Z{\mathbb Z}

\def\D{\mathcal D}
\def\wedgeL{\wedge_L}

\def\veeL{\vee_{\!L}}

\def\fr{\operatorname{frac}}
\let\le\leqslant
\let\ge\geqslant
\def\len{\operatorname{len}}
\def\infs{\inf{\!}_ s}
\def\sups{\sup{\!}_s}
\def\lens{\len{\!}_s}

\def\t{\operatorname{\it t}}
\def\Lmax{\operatorname{L^{max}}}

\def\INF{\t_{\inf}}
\def\SUP{\t_{\sup}}
\def\LEN{\t_{\len}}

\begin{document}

\title[Translation numbers in a Garside group]
{Translation numbers in a Garside group are
rational with uniformly bounded denominators}

\author{Eon-Kyung Lee \and Sang Jin Lee}
\address{Department of Applied Mathematics, Sejong University,
    Seoul, 143-747, Korea}
\email{eonkyung@sejong.ac.kr}
\address{Department of Mathematics, Konkuk University,
    Seoul, 143-701, Korea}
\email{sangjin@konkuk.ac.kr}

\begin{abstract}
It is known that Garside groups are strongly translation discrete.
In this paper, we show that
the translation numbers in a Garside group are rational with uniformly
bounded denominators and can be computed in finite time.
As an application, we give solutions to some group-theoretic problems.

\end{abstract}

\maketitle

\section{Introduction}\label{sec:introduction}

The notion of translation numbers
was first introduced by Gersten and Short~\cite{GS91}.
For a finitely generated group $G$ and
a finite set $X$ of semigroup generators for $G$,
the \emph{translation number} with respect to $X$ of
a non-torsion element $g\in G$ is defined by
$$\t_{G,X}(g)=\lim_{n\to \infty}\frac{|g^n|{}_X}n,$$
where $|\cdot|_X$ denotes the shortest word length
in the alphabet $X$.
If there is no confusion about the group $G$,
we simply write $\t_X(g)$ instead of $\t_{G,X}(g)$.
When $A$ is a set of group generators, $|g|_A$ and $\t_A(g)$ indicate
$|g|_{A\cup A^{-1}}$ and $\t_{A\cup A^{-1}}(g)$, respectively.

The notion of translation numbers is quite a useful one
since it has both algebraic and geometric aspects.
Kapovich~\cite{Kap97} and Conner~\cite{Con00}
suggested the following notions:
a finitely generated group $G$ is said to be
\begin{enumerate}
\item \emph{translation separable} (or \emph{translation proper})
if for some (and hence for any) finite set $X$
of semigroup generators for $G$
the translation numbers of non-torsion elements
are strictly positive;
\item \emph{translation discrete} if it is translation separable
and for some (and hence for any) finite set $X$
of semigroup generators for $G$
the set $\t_X(G)$ has 0 as an isolated point;
\item \emph{strongly translation discrete} if it is translation separable
and for some (and hence for any) finite set $X$
of semigroup generators for $G$
and for any real number $r$
the number of conjugacy classes $[g]=\{h^{-1}gh:h\in G\}$
with $\t_X(g)\le r$ is finite.
(The translation number is constant on
each conjugacy class~\cite{GS91}.)
\end{enumerate}

There are several results on translation numbers in
geometric and combinatorial groups.
Biautomatic groups are translation separable~\cite{GS91}.
C(4)-T(4)-P, C(3)-T(6)-P and C(6)-P small cancelation
groups are strongly translation discrete~\cite{Kap97}.
Word hyperbolic groups are strongly translation discrete,
and moreover, the translation numbers in a word hyperbolic group are rational
with uniformly bounded denominators~\cite{Gro87,BGSS91,Swe95}.

This paper discusses the translation numbers in Garside groups.
The class of Garside groups, first introduced by Dehornoy
and Paris~\cite{DP99},
provides a lattice-theoretic
generalization of braid groups and Artin groups of finite type.
Bestvina~\cite{Bes99} showed that Artin groups of finite
type are translation discrete.
Charney, Meier and Whittlesey~\cite{CMW04}
showed that Garside groups with tame Garside element
are translation discrete.
Recently, Lee~\cite{Lee07} showed that
Garside groups are strongly translation discrete
without any assumption.

The goal of this paper is to establish the following result
which asserts that, as in word hyperbolic groups,
the translation numbers in each Garside group are rational
with uniformly bounded denominators.
(See the next section for the definitions of $\Delta$, $\D$
and $\Vert\cdot\Vert$.)

\par\medskip\noindent\textbf{Main Theorem}
(Theorem~\ref{thm:trans-rational}~(i))\ \
Let $G$ be a Garside group with Garside element $\Delta$
and the set $\D$ of simple elements.
For every element $g$ of $G$, the translation number
$\t_\D(g)$ is rational of the form $p/q$
for some integers $p, \ q$ such that $1\le q \le \Vert\Delta\Vert^2$.
\par\medskip

For elements of Garside groups,
there are integer-valued invariants, $\inf$ and $\sup$,
that directly come from the definition of Garside groups.
Because $|g|_\D$ is either $-\inf(g)$, $\sup(g)$ or $\sup(g)-\inf(g)$,
it is natural to consider the following limits:
$$
\INF(g)=\lim_{n\to\infty}\frac{\inf(g^n)}n\quad\mbox{and}\quad
\SUP(g)=\lim_{n\to\infty}\frac{\sup(g^n)}n.
$$
Then the translation number $\t_\D(g)$ is either
$-\INF(g)$, $\SUP(g)$ or $\SUP(g)-\INF(g)$.
Like translation numbers, $\INF$ and $\SUP$ are conjugacy invariants.
Exploiting the theory of conjugacy classes in Garside groups,
we study the properties of $\INF$ and $\SUP$, and then prove our Main Theorem.

In addition,
we consider the quotient group
$G_\Delta = G/\langle\Delta^{m_0}\rangle$ of a Garside group $G$,
where $\langle\Delta^{m_0}\rangle$ is a subgroup of the center of $G$.
We show that the group $G_\Delta$ is strongly translation discrete
and the translation numbers in it are rational with uniformly bounded
denominators.

As an application, we show that the power (conjugacy) problem
and the proper power (conjugacy) problem
are solvable in Garside groups and
the generalized power (conjugacy) problem
is solvable in braid groups and Artin groups of type $B$.

\section{Garside groups}\label{sec:GarsideGroups}
We briefly review the definition of Garside groups, and
some results necessary for this work.
See~\cite{Gar69,BS72,Thu92,EM94,BKL98,DP99,Pic01b,Deh02,Geb05,LL06b} for details.

\subsection{Garside monoids and groups}

Let $M$ be a monoid.
Let \emph{atoms}  be the elements $a\in M\setminus \{1\}$
such that $a=bc$ implies either $b=1$ or $c=1$.
Let $\Vert a\Vert$ be the supremum of the lengths of all expressions of
$a$ in terms of atoms. The monoid $M$ is said to be \emph{atomic}
if it is generated by its atoms and $\Vert a\Vert<\infty$ for any $a\in M$.
In an atomic monoid $M$, there are partial orders $\le_L$ and $\le_R$:
$a\le_L b$ if $ac=b$ for some $c\in M$;
$a\le_R b$ if $ca=b$ for some $c\in M$.

\begin{defn}
An atomic monoid $M$ is called a \emph{Garside monoid} if
\begin{enumerate}
\item[(i)] $M$ is finitely generated;
\item[(ii)] $M$ is left and right cancellative;
\item[(iii)] $(M,\le_L)$ and $(M,\le_R)$ are lattices;
\item[(iv)] there exists an element $\Delta$, called a
\emph{Garside element}, satisfying the following:\\
(a) for each $a\in M$, $a\le_L\Delta$ if and only if $a\le_R\Delta$;\\
(b) the set $\{a\in M: a \le_L\Delta\}$ generates $M$.
\end{enumerate}
\end{defn}

An element $a$ of $M$ is called a \emph{simple element} if $a\le_L\Delta$.
Let $\D$ denote the set of all simple elements.
Let $\wedgeL$ and $\veeL$ denote the gcd and lcm with respect to $\le_L$.

Garside monoids satisfy Ore's conditions,
and thus embed in their groups of fractions.
A \emph{Garside group} is defined as the group of fractions
of a Garside monoid.
When $M$ is a Garside monoid and $G$ the group of fractions of $M$,
we identify the elements of $M$ and their images in $G$
and call them \emph{positive elements} of $G$.
$M$ is called the \emph{positive monoid} of $G$,
often denoted $G^+$.

Let $\tau\colon G\to G$ be the inner automorphism of $G$ defined by
$\tau(g)=\Delta^{-1} g\Delta$.
It is known that  $\tau(G^+)=G^+$, that is, the positive monoid is
invariant under the conjugation by $\Delta$.

The partial orders $\le_L$ and $\le_R$, and thus the lattice structures
in the positive monoid $G^+$ can be extended
to the Garside group $G$ as follows:
$g\le_L h$ (respectively, $g\le_R h$) for $g,h\in G$
if $gc=h$ (respectively, $cg=h$) for some $c\in G^+$.

For $g\in G$, there are integers $r\le s$ such that
$\Delta^r\le_L g\le_L\Delta^s$.
Hence, the invariants
$\inf(g)=\max\{r\in\Z:\Delta^r\le_L g\}$,
$\sup(g)=\min\{s\in\Z:g\le_L \Delta^s\}$ and
$\len(g)=\sup(g)-\inf(g)$
are well-defined.
For $g\in G$, there is a unique expression
$$
g=\Delta^r s_1\cdots s_k,
$$
called the \emph{normal form} of $g$,
where $s_1,\ldots,s_k\in \D\setminus\{1,\Delta\}$ and
$(s_is_{i+1}\cdots s_k)\wedgeL \Delta=s_i$ for $i=1,\ldots,k$.
In this case, $\inf(g)=r$ and\/ $\sup(g)=r+k$.

For $g\in G$, we denote its conjugacy class $\{ h^{-1}gh : h\in G\}$
by $[g]$.
Define $\infs(g)=\max\{\inf(h):h\in [g]\}$ and
$\sups(g)=\min\{\sup(h):h\in [g]\}$.
The \emph{super summit set} $[g]^S$
and the \emph{stable super summit set} $[g]^{St}$
are subsets of the conjugacy class of $g$
defined as follows:
\begin{eqnarray*}
[g]^S&=&\{h\in [g]:\inf(h)=\infs(g) \ \mbox{ and } \sup(h)=\sups(g)\};\\{}
[g]^{St} &=&\{h\in [g]^S:h^k\in[g^k]^S \ \mbox{ for all positive integers $k$}\}.
\end{eqnarray*}
Both of these sets are finite and nonempty~\cite{EM94,LL06b,BGG06a}.

In the rest of the paper, if it is not specified,
$G$ is assumed to be a Garside group, whose positive monoid is $G^{+}$,
with Garside element $\Delta$ and the set $\D$ of simple elements,
where $\| \Delta\|$ is simply written as $N$.

\subsection{Some results}

For $a\in G^+$, define $\Lmax(a)$ by $\Lmax(a)=\Delta\wedgeL a$.

\begin{lem}\label{lem:Lmax}
For $a, b\in G^{+}$ and\/ $1\le_L s\le_L\Delta$,
\begin{enumerate}
\item[(i)] $\Lmax(ab) = \Lmax(a \Lmax(b))$;
\item[(ii)] $\Lmax(\tau(a)) = \tau(\Lmax(a))$;
\item[(iii)] $s\le_L\Lmax(sa)$.
\end{enumerate}
\end{lem}

\begin{lem}\label{thm:inf-basic}
For $g,h\in G$,
\begin{enumerate}
\item[(i)] $\inf(gh)\ge \inf(g)+\inf(h)$;
\item[(ii)] $\inf(g)-\len(h)\le\inf(h^{-1}gh)\le\inf(g)+\len(h)$.
\end{enumerate}
\end{lem}

\begin{proof}
(i)\ \
Let $g=\Delta^u a$ and $h=\Delta^v b$, where
$u=\inf(g)$, $v=\sup(h)$ and $a,b\in G^+$.
Since $gh=\Delta^u a\Delta^v b=\Delta^{u+v}\tau^v(a)b$,
we have $\Delta^{u+v}\le_L gh$, and hence
$\inf(gh)\ge u+v=\inf(g)+\inf(h)$.

\medskip
(ii)\ \
Let $h=\Delta^u s_1\cdots s_k$ be the normal form of $h$,
and let $g_1 = \Delta^{-u}g\Delta^u$.
Since $\inf(\tau(x))=\inf(x)$ for any $x\in G$,
\begin{equation}\label{eqn:St-1}
\inf(g_1)=\inf(\tau^u(g))=\inf(g).
\end{equation}

If $\len(h)=k=0$, then we are done by (\ref{eqn:St-1}).

Hence, we may assume $\len(h)=k\ge 1$.
Observe that, for any $g_2 \in G$ and $s\in\D\setminus\{1, \Delta\}$,
\begin{eqnarray*}
\inf(s^{-1}g_2 s) & \ge & \inf(s^{-1})+\inf(g_2)+\inf(s) \ = \ \inf(g_2) -1; \\
\inf(g_2) & = &\inf(s(s^{-1}g_2 s)s^{-1})
 \ \ge \ \inf(s)+\inf(s^{-1}g_2 s)+\inf(s^{-1}) \\
 & = & \inf(s^{-1}g_2 s)-1,
\end{eqnarray*}
whence $\inf(g_2 )-1\le\inf(s^{-1}g_2 s)\le\inf(g_2 )+1$.
Therefore,
\begin{equation}\label{eqn:St-2}
\inf(g_1)-k\le\inf(s_k^{-1}\cdots s_1^{-1}g_1 s_1\cdots s_k)\le\inf(g_1 )+k.
\end{equation}
From (\ref{eqn:St-1}) and (\ref{eqn:St-2}), we obtain the desired result.
\end{proof}

\begin{prop}[Theorem 6.1 of \cite{Lee07};
Lemma 3.5 and Proposition 3.6 of~\cite{LL06b}]
\label{thm:inequality}
For $g\in G$ and $n, m\ge 1$,
\begin{enumerate}
\item[(i)]  $n\infs(g)\le \infs(g^n)\le n\infs(g)+n-1$;
\item[(ii)] $n\sups(g)-(n-1)\le\sups(g^n)\le n\sups(g)$;
\item[(iii)] $\infs(g^m)+\infs(g^n)\le\infs(g^{m+n})
\le\infs(g^m)+\infs(g^n)+1$.
\end{enumerate}
\end{prop}

\section{Properties of $\INF$, $\SUP$ and $\LEN$}\label{sec:linf}

\begin{defn}
For an element $g$ of a Garside group $G$, define
$$
\INF(g)=\limsup_{n\to\infty} \frac{\inf(g^n)}n;\quad
\SUP(g)=\liminf_{n\to\infty} \frac{\sup(g^n)}n;\quad
\LEN(g)=\SUP(g)-\INF(g).
$$
\end{defn}

Since $\inf(g)\le \inf(g^n)/n \le \sup(g^n)/n \le\sup(g)$
for all $n\ge 1$,
both $\INF(g)$ and $\SUP(g)$ are finite-valued.
In fact, we shall see in Lemmas~\ref{thm:INF-basic} and~\ref{thm:SUP-basic}
that $\INF(g)=\lim_{n\to\infty}\inf(g^n)/n$
and $\SUP(g)=\lim_{n\to\infty}\sup(g^n)/n$,
hence $\LEN(g)=\lim_{n\to\infty} \len(g^n)/n$.

Note that,  for all $g\in G$, $\sup(g)= -\inf(g^{-1})$,
whence $\SUP(g)=-\INF(g^{-1})$.
Therefore, we may focus only on $\INF$ in order to know about
$\INF$, $\SUP$ and $\LEN$.

We first explore elementary properties of $\INF$.

\begin{lem}\label{thm:INF-basic}
Let $g,h\in G$.
\begin{enumerate}
\item[(i)] $\INF(h^{-1}gh)=\INF(g)$.
\item[(ii)] $\INF(g)=\lim_{n\to\infty} \inf(g^n)/n$.
\item[(iii)] $\infs(g) \le \INF(g) \le \infs(g)+1$.
\item[(iv)] For all $n\ge 1$, $\INF(g^n)=n\INF(g)$.
\end{enumerate}
\end{lem}

\begin{proof}
(i)\ \
By Lemma~\ref{thm:inf-basic}~(ii),
$\inf(g^n)-\len(h)\le\inf(h^{-1}g^nh)\le\inf(g^n)+\len(h)$.
Dividing by $n$ and then taking upper limits, we get $\INF(g)=\INF(h^{-1}gh)$.

\medskip
In view of (i), we may assume that $g$ belongs to its stable
super summit set, and hence
\begin{equation}\label{eqn:St-3}
\infs(g^k)=\inf(g^k) \quad\mbox{for all}\quad k\ge 1.
\end{equation}

\medskip (ii)\ \
From (\ref{eqn:St-3}) and Proposition~\ref{thm:inequality}~(iii),
$\inf(g^m)+\inf(g^n)\le\inf(g^{m+n})\le\inf(g^m)+\inf(g^n)+1$
for all $m,n\ge 1$.
It is well known that if $a_{m+n}\le a_{m} + a_n$ for all $m, n\ge 0$,
where $a_n\ge 0$ for all $n$,
then $\lim_{n\to\infty} (a_n/n)$ exists~\cite[pp.~189]{AL93}.
The result follows when $a_n=\inf(g^n) -n\inf(g)+1$.

\medskip (iii)\ \
From (\ref{eqn:St-3}) and Proposition~\ref{thm:inequality}~(i),
$\infs(g)\le \inf(g^n)/n< \infs(g)+1$ for all $n\ge 1$.

\medskip (iv)\ \
$\INF(g^n)=\lim_{k\to\infty} n(\inf(g^{kn})/kn)=n\INF(g)$.
\end{proof}

In an analogous way, we have the following properties for $\SUP$.
\begin{lem}\label{thm:SUP-basic}
Let $g,h\in G$.
\begin{enumerate}
\item[(i)] $\SUP(h^{-1}gh)=\SUP(g)$.
\item[(ii)] $\SUP(g)=\lim_{n\to\infty} \sup(g^n)/n$.
\item[(iii)] $\sups(g)-1\le\SUP(g)\le\sups(g)$.
\item[(iv)] For all $n\ge 1$, $\SUP(g^n)=n\SUP(g)$.
\end{enumerate}
\end{lem}

The bounds of $\INF(g)$ and $\SUP(g)$
in Lemma~\ref{thm:INF-basic} (iii) and Lemma~\ref{thm:SUP-basic} (iii)
are not sharp.
Those bounds are improved in ~\cite{LL06a}.
Corollary 3.5 in~\cite{LL06a} shows that
$\infs(g) \le \INF(g) \le \infs(g)+1-1/N$ and $\sups(g)-1+1/N\le\SUP(g)\le\sups(g)$.
Example 3.6 in~\cite{LL06a} shows that these bounds are optimal.

\begin{defn}
An element $g$ of a Garside group is said to be
\begin{enumerate}
\item[(i)] \emph{inf-straight} if $\inf(g)=\INF(g)$;
\item[(ii)] \emph{sup-straight} if $\sup(g)=\SUP(g)$.
\end{enumerate}
\end{defn}

The following proposition provides conditions equivalent to
inf-straightness.
The implication from (ii) to (iii)
is benefited from discussions with Ki Hyoung Ko.

\begin{prop}\label{thm:inf-straight-equiv}
For every $g\in G$, the following conditions are equivalent.
\begin{enumerate}
\item[(i)]
$g$ is inf-straight.
\item[(ii)]
$\inf(g^N)=N\inf(g)$.
\item[(iii)]
$\inf(g^k)=k\inf(g)$ for all $k\ge 1$.
\end{enumerate}
\end{prop}

\begin{proof}
(i) $\Rightarrow$ (ii)\ \
By Lemma~\ref{thm:INF-basic} and the definition of $\infs$,
$$N\inf(g)=N\INF(g) =\INF(g^N)\ge\infs(g^N)\ge \inf(g^N)\ge N\inf(g).$$

\medskip (ii) $\Rightarrow$ (iii)\ \
If $\len(g)=0$, the result is trivial.
So, we may assume $\len(g) \ge 1$.

Let $\inf(g)=r$.
Notice that $g^k\Delta^{-kr}$ is a positive element for $k\ge 1$,
since $\inf(g^k)\ge k\inf(g)=kr$.
For $k\ge 1$, let $s_k=\Lmax(g^k\Delta^{-kr})$
and $a_k=s_k^{-1} g^k\Delta^{-kr}$, hence
$$
g^k=s_ka_k\Delta^{kr}\quad\mbox{and}\quad s_k=\Lmax(s_ka_k).
$$
Since $\inf(g^k)=\inf(s_ka_k)+kr=\inf(s_ka_k)+k\inf(g)$ and $s_k=\Lmax(s_ka_k)$,
$\inf(g^k)=k\inf(g)$ if and only if $\inf(s_ka_k)=0$ or, equivalently, if $s_k\ne\Delta$.
Therefore, it suffices to show that $s_k\ne\Delta$ for all $k\ge 1$.
Let $\psi=\tau^{-r}$.
(Then, $\Delta^r h=\psi(h)\Delta^r$ for all $h\in G$.)

\medskip
First, we claim that for $k\ge 1$,
\begin{eqnarray}
s_{k+1}&=&\Lmax(s_1a_1\psi(s_k)),\label{eqn:St-4}\\
s_{k+1}&=&\Lmax(s_ka_k\psi^k(s_1))\label{eqn:St-4a}.
\end{eqnarray}
Observe that for all $k\ge 1$,
$$
s_{k+1}a_{k+1}\Delta^{(k+1)r}=g^{k+1}=g\cdot g^k
=(s_1a_1\Delta^r)(s_ka_k\Delta^{kr})
= s_1a_1\psi(s_ka_k)\Delta^{(k+1)r}.
$$
Therefore,
$$
s_{k+1}a_{k+1}=s_1a_1\psi(s_ka_k)=s_1a_1\psi(s_k)\psi(a_k).
$$
Since $s_k=\Lmax(s_ka_k)$,
we have $\psi(s_k)=\Lmax(\psi(s_k)\psi(a_k))$ by Lemma~\ref{lem:Lmax}~(ii).
Therefore,
\begin{eqnarray*}
s_{k+1}
&=& \Lmax(s_{k+1}a_{k+1})\\
&=& \Lmax(s_1a_1\psi(s_k)\psi(a_k))\\
&=& \Lmax(s_1a_1\Lmax(\psi(s_k)\psi(a_k)))
    \qquad(\mbox{by Lemma~\ref{lem:Lmax}~(i)})\\
&=& \Lmax(s_1a_1\psi(s_k)).
\end{eqnarray*}
Hence, Equation~(\ref{eqn:St-4}) is proved.
Applying the same argument to $g^{k+1}=g^k\cdot g$,
we obtain Equation~(\ref{eqn:St-4a}).

\medskip
By Equation~(\ref{eqn:St-4a}) and Lemma~\ref{lem:Lmax}~(iii),
we have $s_k\le_L s_{k+1}$ for all $k\ge 1$.
Moreover, we know that $s_1 \neq 1$ because $\len(g)\ge 1$.
Therefore the sequence $s_1, s_2, \ldots$ satisfies
\begin{equation}\label{eq:seq}
1\neq s_1\le_L s_2\le_L \cdots \le_L s_N\le_L\cdots\le_L\Delta.
\end{equation}
Since $\inf(g^N)=rN$ by the hypothesis, we have $s_N\ne\Delta$,
hence $\Vert s_N\Vert$ is strictly less than $N=\Vert\Delta\Vert$.
It follows that
$$
1\le \Vert s_1\Vert \le \Vert s_2\Vert \le
 \cdots \le \Vert s_{N-1}\Vert \le \Vert s_N\Vert \le N-1.
$$
By the pigeonhole principle, there exists $m\in \{1,\ldots, N-1\}$
such that $\Vert s_m\Vert =\Vert s_{m+1}\Vert$, hence $s_m=s_{m+1}$.

\medskip

Now, we claim that for each $j\ge 1$,
\begin{equation}\label{eqn:induct}
s_m=s_{m+j}.
\end{equation}
Using induction, it suffices to show that
if $s_k=s_{k+1}$ for some $k\ge 1$, then $s_{k+1}=s_{k+2}$.
Suppose $s_k=s_{k+1}$.
By Equation~(\ref{eqn:St-4}),
$$
s_{k+1}
= \Lmax(s_1a_1\psi(s_k))
=\Lmax(s_1a_1\psi(s_{k+1}))
=s_{k+2}.
$$
Therefore, Equation~(\ref{eqn:induct}) is proved.

\medskip
Since $s_N\ne\Delta$, one has $s_k\ne\Delta$ for $k\le N$
by Equation~(\ref{eq:seq}).
In particular, $s_m\ne\Delta$ because $m<N$.
By Equation~(\ref{eqn:induct}),
we can conclude that $s_k\ne\Delta$ for all $k\ge m$.
Therefore, $s_k\ne\Delta$ for all $k\ge 1$
and we are done.

\medskip
(iii) $\Rightarrow$ (i)\ \
It is obvious by the definition of $\INF$.
\end{proof}

In an analogous way, we have conditions equivalent to sup-straightness.
\begin{prop}\label{thm:sup-straight-equiv}
For every $g\in G$, the following conditions are equivalent.
\begin{enumerate}
\item[(i)] $g$ is sup-straight.
\item[(ii)]
$\sup(g^N)=N\sup(g)$.
\item[(iii)]
$\sup(g^k)=k\sup(g)$ for all $k\ge 1$.
\end{enumerate}
\end{prop}

Using Proposition~\ref{thm:inf-straight-equiv}~(ii), it is easy
to decide whether or not a given element is inf-straight.
We remark that the value $N$ in that statement is the smallest one playing
such a role.
We illustrate this with an example.

\begin{exmp}\label{ex:1}
For $N\ge M\ge 2$, let
$$G=\langle x,y\mid x^N=y^M\rangle.$$
It is a Garside group with Garside element
$\Delta=x^N=y^M$~\cite[Example 4]{DP99},
and $\Vert\Delta\Vert=N$.
Note that $\inf(x^k)=\lfloor k/N\rfloor$.
Hence, $\inf(x^k)=0=k\inf(x)$ for all $1\le k\le N-1$ but
$x$ is not inf-straight because
$\INF(x)=1/N$ is not equal to $\inf(x)=0$.
\end{exmp}

\begin{lem}\label{thm:inf-straight}
If $g\in G$ is inf-straight, then $g^n$ is also inf-straight
and $\infs(g^n)=n\inf(g)$ for all $n\ge 1$.
\end{lem}

\begin{proof}
Let $g\in G$ be inf-straight, and let $n\ge 1$.
By Lemma~\ref{thm:INF-basic} and Proposition~\ref{thm:inf-straight-equiv},
$$
\INF(g^n)=n\INF(g)=n\inf(g)=\inf(g^n)\le\infs(g^n)\le\INF(g^n).
$$
In particular, $\infs(g^n)=n\inf(g)$.
Since $\INF(g^n)=\inf(g^n)$, $g^n$ is inf-straight.
\end{proof}

\begin{lem}\label{thm:sup-straight}
If $g\in G$ is sup-straight, then $g^n$ is also sup-straight and
$\sups(g^n)=n\sup(g)$ for all positive integers $n$.
\end{lem}

The stable super summit set is nonempty~\cite{LL06b}.
Using this fact together with the properties of inf-straight elements
and sup-straight elements,
we characterize conjugacy classes of inf-straight elements and sup-straight elements.

\begin{cor}\label{thm:inf-straight-conj}
For every $g\in G$, the following conditions are equivalent.
\begin{enumerate}
\item[(i)]
$g$ is conjugate to an inf-straight element.
\item[(ii)]
$\infs(g)=\INF(g)$.
\item[(iii)]
$\infs(g^N)=N\infs(g)$.
\item[(iv)]
$\infs(g^k)=k\infs(g)$ for all $k\ge 1$.
\item[(v)]
For all $h\in[g]^S$, $h$ is inf-straight.
\end{enumerate}
\end{cor}

\begin{proof}
(i) $\Rightarrow$ (ii)\ \
Suppose that $h\in [g]$ is inf-straight.
By Lemmas~\ref{thm:INF-basic} and~\ref{thm:inf-straight},
$$\INF(g)=\INF(h)=\inf(h)=\infs(h)=\infs(g).$$

\medskip
(ii) $\Leftrightarrow$ (iii) $\Leftrightarrow$ (iv)\ \
Because $\infs$ and $\INF$ are conjugacy invariants,
we may assume that $g$ belongs to its stable super summit set.
Then $\infs(g^k)=\inf(g^k)$ for all $k\ge 1$, hence
the statements (ii), (iii) and (iv)
are exactly the same as those in Proposition~\ref{thm:inf-straight-equiv}.

\medskip
(iv) $\Rightarrow$ (v)\ \
Let $h\in[g]^S$.
For all $k\ge 1$,
$$
k\inf(h)\le \inf(h^k)\le\infs(g^k)=k\infs(g)=k\inf(h).
$$
Therefore $\inf(h^k)=k\inf(h)$ for all $k\ge 1$, hence $h$ is
inf-straight by Proposition~\ref{thm:inf-straight-equiv}.

\medskip
(v) $\Rightarrow$ (i)\ \
It is obvious.
\end{proof}

\begin{cor}\label{thm:sup-straight-conj}
For every $g\in G$, the following conditions are equivalent.
\begin{enumerate}
\item[(i)]
$g$ is conjugate to a sup-straight element.
\item[(ii)]
$\sups(g)=\SUP(g)$.
\item[(iii)]
$\sups(g^N)=N\sups(g)$.
\item[(iv)]
$\sups(g^k)=k\sups(g)$ for all $k\ge 1$.
\item[(v)]
For all $h\in[g]^S$, $h$ is sup-straight.
\end{enumerate}
\end{cor}

For a real number $x$, let $\fr(x)$ denote the fractional part of $x$, that is,
$\fr(x)=x-\lfloor x\rfloor$.

\begin{lem}\label{thm:Ko}
For every $g\in G$, $\fr\left(\INF(g)\right)\not\in (0, 1/N)$.
\end{lem}

\begin{proof}
Assume that $0<\fr(\INF(g))< 1/N$ for some $g\in G$.
Since $\infs(g)\le\INF(g)\le\infs(g)+1$ by Lemma~\ref{thm:INF-basic}
and $\infs(g)$ is integer-valued, we have
\begin{equation}\label{eqn:St-5}
\infs(g) < \INF(g)<\infs(g)+\frac1N.
\end{equation}
Note that $\infs(g^N)\ge N\infs(g)$
by Proposition~\ref{thm:inequality}~(i).
If $\infs(g^N)\ge N\infs(g)+1$, then by Lemma~\ref{thm:INF-basic}
$$
\INF(g) = \frac{\INF(g^N)}N
\ge \frac{\infs(g^N)}N
\ge \frac{N\infs(g)+1}N
=\infs(g)+\frac1N,
$$
which contradicts (\ref{eqn:St-5}).
Therefore $\infs(g^N)= N\infs(g)$ and it follows
by Corollary~\ref{thm:inf-straight-conj} that $\INF(g)=\infs(g)$.
This contradicts (\ref{eqn:St-5}).
\end{proof}

We are now ready to show the main result of this section.

\begin{thm}\label{thm:INF-main}
Let $G$ be a Garside group with Garside element $\Delta$, and
let $N=\|\Delta\|$.
\begin{enumerate}
\item[(i)]
For every element $g$ of $G$,
$\INF(g)$ and $\SUP(g)$ are rational of the form $p/q$
for some integers $p, \ q$ such that $1\le q\le N$.
\item[(ii)]
For every element $g$ of $G$,
$\LEN(g)$ is rational of the form $p/q$
for some integers $p, \ q$ such that $1\le q\le N^2$.
\item[(iii)]
There is a finite-time algorithm that, given an element $g$ of\/ $G$,
computes $\INF(g)$, $\SUP(g)$ and $\LEN(g)$.
\end{enumerate}
\end{thm}

\begin{proof}
(i)\ \
Put $\INF(g)=x$.
Assume that $x$ is irrational.
Then, the set $\{\fr(kx) : k\ge 1\}$ is a dense subset of
the unit interval $[0,1]$.
Therefore there exists $k\ge 1$ such that
$0<\fr(kx) < 1/N$.
Since $kx=k\INF(g)=\INF(g^k)$, it contradicts Lemma~\ref{thm:Ko}.
Consequently, $x$ is rational.

If $x=0$, we are done.
Thus, we may assume $x\neq 0$.
Let $x=p/q$ for relatively prime integers $p,q$ with $q\ge1$.
Assume $q>N$.
There exist integers $k\ge 1$ and $m$ such that $kp=mq+1$.
Hence, $\INF(g^k)=k\INF(g)=kp/q=m+1/q \in (m, m+1/N)$, which
contradicts Lemma~\ref{thm:Ko}.
Consequently, $1\le q\le N$.

The assertion for $\SUP(g)$ in the statement is true
because $\SUP(g)=-\INF(g^{-1})$.

\medskip
(ii) follows from (i) because $\LEN(g)=\SUP(g)-\INF(g)$.

\medskip (iii)\ \
If $N=1$, then $G$ is the infinite cyclic group generated
by $\Delta$, and thus $\INF(g)=\inf(g)$ for all $g\in G$.
So, we may assume $N\ge 2$.
To compute $\INF(g)$, we will make an approximation of $\INF(g)$
using $\infs(g^n)/n$ for sufficiently large $n$
(in fact, $n=N^2$ is enough),
and then determine $\INF(g)$ exactly by using (i).

Let $T$ be the set of all rational numbers whose denominators
are less than or eqal to $N$, that is,
$T=\{p/q: p,q\in\Z,\ 1\le q\le N\}$.
Then $\INF(g)\in T$ by (i).
Note that in any closed interval $[\alpha,\alpha+\epsilon]$ on the real line
with length $\epsilon\le 1/N^2$,
there is at most one element of $T$.

Choose any $n\ge N^2$.
Applying Lemma~\ref{thm:INF-basic}~(iii) to $g^n$,
we know $\infs(g^n)\le\INF(g^n)\le\infs(g^n)+1$.
Since $\INF(g^n)=n\INF(g)$ by Lemma~\ref{thm:INF-basic}~(iv),
\begin{equation}\label{eqn:INF-est1}
\frac{\infs(g^n)}n \le \INF(g) \le \frac{\infs(g^n)}n+\frac1n.
\end{equation}
Therefore, $\INF(g)$ belongs to the closed interval
$[\infs(g^n)/n,\ \infs(g^n)/n+1/n]$.
Since the length of this interval is $\le 1/N^2$,
$\INF(g)$ is the unique element of $T$ in this interval,
and hence it can be found in finite time.
Because $\SUP(g)=-\INF(g^{-1})$ and $\LEN(g)=\SUP(g)-\INF(g)$
for all $g\in G$, there is a finite-time algorithm
for $\INF(\cdot)$, $\SUP(\cdot)$ and $\LEN(\cdot)$.
\end{proof}

We close this section with examples concerning the upper bounds
of the denominators of $\INF$, $\SUP$ and $\LEN$
given in Theorem~\ref{thm:INF-main}.
Recall the Garside group
$$
G=\langle x,y\mid x^N=y^M\rangle\qquad (N\ge M\ge 2)
$$
considered in Example~\ref{ex:1}.
For $k\ge 1$, $\inf(x^k)=\lfloor k/N\rfloor$ and
$\sup(x^k)=\lceil k/N\rceil$, whence $\INF(x)=\SUP(x)=1/N$.
Therefore, the upper bound $N$ of the denominators
of $\INF$ and $\SUP$ is optimal.

The following example shows that the upper bound $N^2$
of the denominator of $\LEN(g)$ is asymptotically tight.
That is, the optimal upper bound of the denominator equals $\Theta(N^2)$.

\begin{exmp}\label{ex:tinf}
Consider the group
$$
G=H\times H,\quad
\mbox{where $H=\langle x,y\mid x^p=y^q\rangle$ with $q=p+1\ge 3$.}
$$
As in Example~\ref{ex:1},
$H$ is a Garside group with Garside element $\Delta_H =x^p=y^q$.
Because $G$ is a cartesian product of $H$,
it is a Garside group with Garside element
$\Delta = (\Delta_H, \Delta_H)$ by~\cite[Theorem 4.1]{Lee07}.
It is clear that $N=\Vert\Delta\Vert= 2q$.
Let $g=(x,y)\in G$.
For all $n\ge 1$, $g^n=(x^n,y^n)$,
so by~\cite[Lemma 4.4]{Lee07}
\begin{eqnarray*}
\inf(g^n)
&=&\min\{\inf(x^n),\inf(y^n)\}
 =\min\{\lfloor n/p\rfloor,\lfloor n/q\rfloor\}
 =\lfloor n/q\rfloor;\\
\sup(g^n)
&=&\max\{\sup(x^n),\sup(y^n)\}
 = \max\{\lceil n/p\rceil,\lceil n/q\rceil\}
 = \lceil n/p\rceil.
\end{eqnarray*}
Therefore,
$$
\INF(g)=\frac1q,\quad
\SUP(g)=\frac1p\quad\mbox{and}\quad
\LEN(g)=\frac1p-\frac1q=\frac{q-p}{pq}=\frac1{\frac N2(\frac N2-1)}.
$$
Since $N$ is even, the denominator of $\LEN(g)$ is
$\frac N2(\frac N2-1)=\Theta(N^2)$.
\end{exmp}

\section{Translation numbers}\label{sec:TransNum}

It is well known (see~\cite{Cha95} for example) that,
for an element $g\in G$, the shortest word length $|g|_\D$
can be expressed in terms of $\inf(g)$ or $\sup(g)$ as follows:
(i) if $\inf(g)\ge 0$, then $|g|_{\D}=\sup(g)$;
(ii) if $\sup(g)\le 0$, then $|g|_{\D}=-\inf(g)$;
(iii) if $\inf(g)< 0<\sup(g)$, then $|g|_{\D}=\sup(g)-\inf(g)=\len(g)$.

\begin{lem}\label{thm:trans-INF}
Let $G$ be a Garside group with the set $\D$ of simple elements.
Let $g\in G$.
\begin{itemize}
\item[(i)] If\/ $\infs(g)\ge 0$, then $\t_\D(g)=\SUP(g)$.
\item[(ii)] If\/ $\sups(g)\le 0$, then $\t_\D(g)=-\INF(g)$.
\item[(iii)] If\/ $\infs(g)< 0<\sups(g)$, then $\t_\D(g)=\LEN(g)$.
\end{itemize}
\end{lem}

\begin{proof}
Note that $\infs$, $\sups$, $\t_\D$, $\INF$, $\SUP$ and $\LEN$
are all conjugacy invariants.
Hence, we may assume $g\in[g]^{St}$,
that is, $\infs(g^n)=\inf(g^n)$ and $\sups(g^n)=\sup(g^n)$
for all $n\ge 1$.

If $\infs(g)\ge 0$, then $\inf(g^n)\ge n\inf(g)\ge 0$ for all $n\ge 1$,
whence $|g^n|_\D =\sup(g^n)$.
Consequently, $\t_\D(g)=\SUP(g)$.

If $\sups(g)\le 0$, then $\sup(g^n)\le n\sup(g)\le 0$ for all $n\ge 1$,
whence $|g^n|_\D =-\inf(g^n)$.
Consequently, $\t_\D(g)=-\INF(g)$.

If $\infs(g)<0<\sups(g)$, then
$\inf(g^n) \le n\inf(g) + (n-1) < 0 < n\sup(g) -(n-1) \le \sup(g^n)$
for all $n\ge 1$ by Proposition~\ref{thm:inequality},
whence $|g^n|_\D=\len(g^n)$.
Consequently, $\t_\D(g)=\LEN(g)$.
\end{proof}

\begin{thm}\label{thm:trans-rational}
Let $G$ be a Garside group with Garside element $\Delta$
and the set $\D$ of simple elements. Let $N=\Vert\Delta\Vert$.
\begin{enumerate}
\item[(i)]
The translation numbers in $G$ are
rational of the form $p/q$
for some integers $p, \ q$ such that $1\le q\le N^2$.
\item[(ii)]
If\/ $g$ is a non-identity element of\/ $G$, then $\t_\D(g)\ge 1/N$.
\item[(iii)]
There is a finite-time algorithm that, given an element of\/ $G$,
computes its translation number.
\end{enumerate}
\end{thm}

\begin{proof}
(i) and (iii) immediately follow from Theorem~\ref{thm:INF-main}
and Lemma~\ref{thm:trans-INF}. Let us prove (ii).

Let $g$ be a non-identity element of $G$.
Because Garside groups are torsion-free~\cite{Deh98} and
strongly translation discrete (and hence translation separable)~\cite{Lee07},
$\t_\D(g) >0$.

If $\infs(g)\ge 0$, then $\t_\D(g)=\SUP(g)\ge 1/N$ by
Lemma~\ref{thm:trans-INF} and Theorem~\ref{thm:INF-main}.

If $\sups(g)\le 0$, then $\t_\D(g)=-\INF(g)=|\INF(g)|\ge1/N$ by
Lemma~\ref{thm:trans-INF} and Theorem~\ref{thm:INF-main}.

If $\infs(g)<0<\sups(g)$, then
$\infs(g^n) < 0 < \sups(g^n)$ for all $n\ge 1$
by Proposition~\ref{thm:inequality},
whence $\INF(g)\le 0\le\SUP(g)$.
Therefore, $\t_\D(g)=\SUP(g)-\INF(g)=|\SUP(g)|+|\INF(g)|$
by Lemma~\ref{thm:trans-INF}.
Since $\t_\D(g)> 0$, either $|\SUP(g)|$ or $|\INF(g)|$ is greater than
or equal to $1/N$ by Theorem~\ref{thm:INF-main}.
Consequently, $\t_\D(g)\ge 1/N$.
\end{proof}

In the rest of this section, we study the translation numbers
in the quotient group
$G_\Delta=G/\langle\Delta^{m_0}\rangle$,
where $G$ is a Garside group and $m_0$ is the smallest power of $\Delta$
such that $\Delta^{m_0}$ is central in $G$.
In some cases, this quotient group is more natural and convenient to deal with.
If $G$ is the $n$-string braid group
and $\Delta$ is the usual Garside element,
then $G_\Delta$ is the mapping class group of the $n$-punctured disk.
Both Bestvina~\cite{Bes99} and Charney, Meier and Whittlesey~\cite{CMW04}
studied the cocompact action of this quotient group,
rather than the Garside group itself, on a finite dimensional
complex constructed using the normal form.

For an element $g$ of a Garside group $G$,
let $\bar g$ denote the image of $g$ under
the canonical epimorphism $G\to G_\Delta$.
Let $\bar\D=\{\bar s:s\in\D\}$ and
let $\t_{\bar\D}(\bar g)$
denote the translation number
of $\bar g\in G_\Delta$ with respect to $\bar\D$.

\begin{lem}\label{thm:trans-G_Delta}
For every $g\in G$,
$\t_{\bar\D}(\bar g)=\LEN(g)$.
\end{lem}

\begin{proof}
It is easy to see that for any $h\in G$,
$\len(h)\le|\bar h|_{\bar\D}<\len(h)+m_0$.
Substituting $g^n$ for $h$,
we obtain
$\len(g^n)\le|(\bar g)^n|_{\bar\D}<\len(g^n)+m_0$,
from which the assertion in the statement follows.
\end{proof}

\begin{lem}\label{thm:LEN-basic}
For every $g\in G$, $\lens(g)-2\le\LEN(g)\le\lens(g)$.
\end{lem}

\begin{proof}
It immediately follows from Lemmas~\ref{thm:INF-basic}~(iii)
and~\ref{thm:SUP-basic}~(iii).
\end{proof}

\begin{thm}\label{thm:G_Delta_LB}
Let $G$ be a Garside group with Garside element $\Delta$.
Then $G_\Delta$ is strongly translation discrete.
Moreover, the translation numbers in $G_\Delta$ are rational of the form
$p/q$ for some integers $p, \ q$ such that $1\le q\le N^2$,
where $N=\Vert\Delta\Vert$.
\end{thm}

\begin{proof}
The translation numbers in $G_\Delta$ are rational of the form
as in the statement by Lemma~\ref{thm:trans-G_Delta}
and Theorem~\ref{thm:INF-main}~(ii).

We first show that $G_\Delta$ is translation separable.
Suppose that $\t_{\bar \D}(\bar g)=\LEN(g)=0$ for $g\in G$.
(We shall see that $\bar g$ has finite order.)
Because both the translation numbers and the orders
of elements are conjugacy invariants in $G_\Delta$,
we may assume that $g$ belongs to its
stable super summit set.
That is, for all $n\ge 1$, $\inf(g^n)=\infs(g^n)$
and $\sup(g^n)=\sups(g^n)$, whence $\len(g^n)=\lens(g^n)$.
By Lemmas~\ref{thm:LEN-basic},~\ref{thm:SUP-basic}~(iv)
and~\ref{thm:INF-basic}~(iv),
$$
\len(g^n)\le\LEN(g^n)+2=n\LEN(g)+2=2\qquad\mbox{for all $n\ge 1$.}
$$
Since $\len(g^{-n})=\len(g^n)\le 2$ for all $n\ge 1$,
$$
\{\bar g^n:n\in\Z\}\subset\{\bar h: 0\le \inf(h)<m_0, \ \len(h)\le 2\}.
$$
Because the right hand side is a finite set,
the cyclic group generated by $\bar g$ is a finite subgroup of $G_\Delta$.
Hence, $\bar g$ has finite order in $G_\Delta$.

\medskip
Now we show that $G_\Delta$ is strongly translation discrete.
For a real number $r$, let
\begin{eqnarray*}
A_r&=&\{\bar g\in G_\Delta: \t_{\bar\D}(\bar g)\le r\};\\
C_r&=&\{h\in G:0\le \inf(h) < m_0,\ \len(h)\le r+2\}.
\end{eqnarray*}
Because $G_\Delta$ is translation separable,
it suffices to show that for any real number $r$, there are only finitely
many conjugacy classes in $A_r$.
Because $C_r$ is a finite set, it suffices to show that
for each $\bar g\in A_r$, there exists $h\in C_r$
such that $\bar h$ is conjugate to $\bar g$ in $G_\Delta$.

Let $\bar g\in A_r$ for $g\in G$.
By Lemmas~\ref{thm:LEN-basic} and~\ref{thm:trans-G_Delta}
$$
\lens(g)\le\LEN(g)+2=\t_{\bar\D}(\bar g)+2\le r+2.
$$
Choose any element $h_0$ in the super summit set of $g$.
Let $k$ be such that $0\le \inf(h_0)-km_0< m_0$,
and then let $h=\Delta^{-km_0}h_0$.
Since
$$
0\le\inf(h) \ (=\inf(h_0)-km_0)< m_0\quad\mbox{and}\quad
\len(h) \ (=\len(h_0)=\lens(g))\le r+2,
$$
$h\in C_r$.
By the construction,
$\bar h$ is conjugate to $\bar g$.
\end{proof}

We remark that, for Theorem~\ref{thm:trans-rational},
the upper bound $N^2$ on the denominators of
translation numbers is asymptotically tight, and the lower bound
$1/N$ on the translation numbers of non-identity elements is optimal.
For Theorem~\ref{thm:trans-rational}~(ii), recall the Garside group
$G$ and the element $x$ given in Example~\ref{ex:1}.
In this case, $\t_\D(x)=\SUP(x)=1/N$.
For Theorem~\ref{thm:trans-rational}~(i), recall the Garside group
$G$ given in Example~\ref{ex:tinf}.
Let $g=(x^{-1}, y)$.
Then $\t_\D(g)=\LEN(g)=(p+q)/(pq)$.
Since $p+q$ and $pq$ are relatively prime,
the denominator of $\t_\D(g)$ is equal to $pq=(N/2)(N/2-1)=\Theta(N^2)$.
This example also shows the asymptotic tightness of the upper bound $N^2$
given in Theorem~\ref{thm:G_Delta_LB}.

\section{Some group-theoretic problems}

In this section, we apply our results on the translation numbers
to solve some group-theoretic problems in Garside groups.
Lipschutz and Miller~\cite{LM71} considered
the following fundamental problems in groups.
\begin{itemize}
\item The order problem:
given $g\in G$, find $n\ge 1$ such that $g^n=1$.
\item The root problem:
given $g\in G$ and $n\ge 1$, find $h\in G$ such that $h^n=g$.
\item The power problem:
given $g,h\in G$, find $n\in\Z$ such that $h^n=g$.
\item The proper power problem:
given $g\in G$, find $h\in G$ and $n\ge 2$ such that $h^n=g$.
\item The generalized power problem:
given $g,h\in G$, find $n, m\in\Z\setminus\{ 0\}$ such that $g^n=h^m$.
\item The intersection problem for cyclic subgroups:
given $g,h\in G$, find $n, m\in\Z$ such that $g^n=h^m\ne 1$.
\end{itemize}

Because Garside groups are torsion-free~\cite{Deh98}
and have solvable word problem,
the order problem is trivial and
the intersection problem for cyclic subgroups
is equivalent to the generalized power problem.

In addition to the above problems, we consider their conjugacy versions.
\begin{itemize}
\item The root conjugacy problem:
given $g\in G$ and $n\ge 1$, find $h\in G$ such that $h^n$ is conjugate to $g$.
\item The power conjugacy problem:
given $g,h\in G$, find $n\in\Z$ such that $h^n$ is conjugate to $g$.
\item The proper power conjugacy problem:
given $g\in G$, find $h\in G$ and $n\ge2$ such that
$h^n$ is conjugate to $g$.
\item The generalized power conjugacy problem:
given $g,h\in G$, find $n, m\in\Z\setminus\{0\}$ such that
$g^n$ is conjugate to $h^m$.
\end{itemize}

Because the conjugacy problem is solvable in Garside groups,
the root problem is equivalent to the root conjugacy problem.
That is, the root problem for $(g,n)$ is solvable if and only
if so is the root conjugacy problem for $(g,n)$.
Moreover, it is easy to draw a solution to one problem from
a solution to the other problem.
(For example, if $h^n$ is conjugate to $g$, then
we can find $x\in G$ such that $x^{-1}h^n x=(x^{-1}hx)^n=g$,
hence $x^{-1}hx$ is a solution to the root problem for $(g,n)$.)
Similarly, the proper power problem is equivalent to
its conjugacy version.

The root problem is solvable in braid groups by Sty\v shnev~\cite{Sty78}
and in Garside groups by Sibert~\cite{Sib02} (under a mild assumption)
and by Lee~\cite{Lee07} (without any assumption).
Therefore, the root conjugacy problem is also solvable in Garside groups.

We now apply our results on translation numbers
to solve the power, proper power
and generalized power problems in Garside groups
together with their conjugacy versions.

It is easy to see that for elements $g$ and $h$ of a Garside group,
\begin{itemize}
\item[(i)] $\t_\D(g)=\t_\D(g^{-1})$;
\item[(ii)] $\t_\D(g^n)=|n|\t_\D(g)$ for $n\in\Z$;
\item[(iii)] If $h\ne 1$ (hence $\t_\D(h)\ne 0$ because Garside groups
are torsion free and translation separable)
and $h^n$ is conjugate to $g$ for some $n\in\Z$, then
$|n|=\t_\D(g)/\t_\D(h)$.
\end{itemize}

\begin{thm}
The power problem and the power conjugacy problem are solvable in Garside groups.
\end{thm}

\begin{proof}
It is a direct consequence of Theorem~\ref{thm:trans-rational}.
We prove only for the power conjugacy problem.
The power problem can be solved in almost the same way.

Let $G$ be a Garside group.
Suppose we are given $g,h\in G$ and want to find $n\in\Z$ such that
$h^n$ is conjugate to $g$.
We may assume that $g, h\ne 1$, otherwise the problem is trivial.
Let $m=\t_\D(g)/\t_\D(h)$, which can be computed in finite time
by Theorem~\ref{thm:trans-rational}~(iii).
Then, $h^n$ is conjugate to $g$ for some $n\in \Z\setminus\{0\}$
if and only if $m$ is a positive integer and
$h^m$ is conjugate to either $g$ or $g^{-1}$.
\end{proof}

\begin{thm}
The proper power problem and the proper power conjugacy problem
are solvable in Garside groups.
\end{thm}

\begin{proof}
It is a direct consequence of the fact
that Garside groups are translation discrete
and have solvable root problem.
Because the two problems are equivalent in Garside groups as we observed before,
we prove only for the proper power conjugacy problem.

Let $G$ be a Garside group.
Suppose we are given $g\in G$ and want to find $h\in G$ and $n\ge2$
such that $h^n$ is conjugate to $g$.
We may assume that $g\neq 1$, otherwise the problem is trivial.
If such $h$ and $n$ exist, then $h\ne 1$ and so
$n=\t_\D(g)/\t_\D(h)\le N\t_\D(g)$ by Theorem~\ref{thm:trans-rational}~(ii).
Therefore, solving the proper power conjugacy problem for $g$
can be reduced to solving the root conjugacy problems
for $g$ and $n$ with $2\le n\le N\t_\D(g)$.
Since the root conjugacy problem is solvable in $G$,
the proper power conjugacy problem is solvable in $G$.
\end{proof}

We say that a group $G$ has the \emph{unique root property}
if for every $g\in G$ and $n\ge 1$, there exists a unique $n$-th root of $g$,
that is, if $g^n=h^n$ for $g,h\in G$ and $n\ge 1$, then $g=h$.

\begin{thm}
If a Garside group $G$ has a finite index subgroup $G_0$
that has the unique root property, then
the generalized power problem and the generalized power conjugacy problem
are solvable in $G$.
\end{thm}

\begin{proof}
Let $r$ be a positive integer such that $g^r\in G_0$ for all $g\in G$.
(It suffices to take $r=[G:G_0]!$, where $[G:G_0]$ denotes
the index of $G_0$ in $G$.)
Suppose that we are given two elements $g,h\in G\setminus\{1\}$.
(If either $g$ or $h$ is equal to 1, the problem is trivial.)
We want to know whether or not $g^n$ is equal/conjugate to $h^m$
for some non-zero integers $n, m$.
By Theorem~\ref{thm:trans-rational},
both $\t_\D(g)$ and $\t_\D(h)$ are positive rational numbers,
whence we can find positive integers $p,q$ such that $p\t_\D(g)=q\t_\D(h)$
in finite time.

\medskip\noindent
\textbf{Claim.}\ \
(i) $g^n=h^m$ for some $n,m\in\Z\setminus\{ 0\}$
if and only if $g^{pr}$ is equal to either $h^{qr}$ or $h^{-qr}$.

(ii) $g^n$ is conjugate to $h^m$ for some $n,m\in\Z\setminus\{ 0\}$
if and only if $g^{pr}$ is conjugate to either $h^{qr}$ or $h^{-qr}$.

\medskip\noindent\textit{Proof of Claim.}\ \
(i)\ \
If $g^{pr}$ is equal to $h^{qr}$ or $h^{-qr}$,
we may take $(pr,qr)$ or $(pr,-qr)$ as $(n,m)$.

Conversely, suppose $g^n=h^m$ for some $n, m \in\Z\setminus\{ 0\}$.
Recall that $p$, $q$ and $r$ are positive integers such that
$p\t_\D(g)=q\t_\D(h)$ and $g^r, h^r \in G_0$.
Since $(g^n)^{pqr}=(h^m)^{pqr}$, we have
\begin{eqnarray}
 & (g^{pr})^{|qn|}=(h^{qr})^{|pm|}\quad\mbox{or}\quad
    (g^{pr})^{|qn|}=(h^{-qr})^{|pm|}; \label{eqn:St-7} \\
 & |qn|\t_\D(g^{pr}) = |pm|\t_\D(h^{qr}). \label{eqn:St-6}
\end{eqnarray}

Because $\t_\D(g^{pr})=pr\t_\D(g)=qr\t_\D(h)=\t_\D(h^{qr})$
and this is strictly positive,
$|qn| = |pm|$ follows from~(\ref{eqn:St-6}).
In addition, $g^{pr}$ and $h^{\pm qr}$ in~(\ref{eqn:St-7}) belong to $G_0$
which has the unique root property.
Consequently,
$g^{pr}$ is equal to either $h^{qr}$ or $h^{-qr}$.

\medskip (ii)\ \
If $g^{pr}$ is conjugate to $h^{qr}$ or $h^{-qr}$,
we may take $(pr,qr)$ or $(pr,-qr)$ as $(n,m)$.

Conversely, suppose $g^n$ is conjugate to $h^m$ for some
$n, m\in\Z\setminus\{ 0\}$.
Then $g^n=x^{-1}h^mx=(x^{-1}hx)^m$ for some $x\in G$.
By (i), $g^{pr}$ is equal to either $(x^{-1}hx)^{qr}=x^{-1}h^{qr}x$ or
$(x^{-1}hx)^{-qr}=x^{-1}h^{-qr}x$.
Consequently, $g^{pr}$ is conjugate to either $h^{qr}$ or $h^{-qr}$.
\hfill \textit{End of Proof of Claim}

\medskip
The generalized power problem
and the generalized power conjugacy problem for $(g,h)$ is solvable
by the above claim, because we can find $p,q$ in a finite number of steps
and the word problem and the conjugacy problem are solvable in Garside groups.
\end{proof}

Observe that pure braid groups have the unique root property:
they are biorderable by Kim and Rolfsen~\cite{KR03}
and biorderable groups have the unique root property.
It is well known that each of the braid groups (a.k.a.{} Artin groups of type $A$)
and Artin groups of type~$B$ contains a pure braid group as a finite index
subgroup.

\begin{cor}
The generalized power problem and the generalized power conjugacy
problem are solvable in braid groups and
Artin groups of type $B$.
\end{cor}

\section*{Acknowledgements}

The authors cordially thank Ki Hyoung Ko,
Won Taek Song, Dong Han Kim and an anonymous
referee for helpful comments on the paper.
Especially, the implication from (ii) to (iii)
in Proposition~\ref{thm:inf-straight-equiv}
is benefited from discussions with Ki Hyoung Ko.
This paper was supported by Konkuk University in 2006.


\begin{thebibliography}{MMMl}

\bibitem{AL93}
L. Alseda, J. Llibre,
Combinatorial dynamics and entropy in dimension one,
World Scientific Publishing Company, 1993.


\bibitem{BGSS91}
G. Baumslag, S. Gersten, M. Shapiro, H. Short,
Automatic groups and amalgams,
J. Pure Appl. Algebra 76 (1991)
229--316.

\bibitem{Bes99}
M.~Bestvina,
Non-positively curved aspects of Artin groups of finite type,
Geometry and Topology 3 (1999) 269--302.

\bibitem{BGG06a}
J.~Birman, V.~Gebhardt, J.~Gonzalez-Meneses,
Conjugacy in Garside groups I: Cyclings, powers, and rigidity,
arXiv:math.GT/0605230.

\bibitem{BKL98}
J.~Birman, K.H.~Ko, S.J.~Lee,
A new approach to the word and conjugacy problems in the braid groups,
Adv. Math. 139 (1998) 322--353.

\bibitem{BS72}
E.{} Brieskorn, K.{} Saito,
Artin-Gruppen und Coxeter-Gruppen,
Invent. Math. 17 (1972) 245--271.

\bibitem{Cha95}
R. Charney,
Geodesic automation and growth functions for Artin groups of finite type,
Math. Ann. 301 (1995) 307-324.

\bibitem{CMW04}
R. Charney, J. Meier, K. Whittlesey,
Bestvina's normal form complex and the homology of Garside groups,
Geom. Dedicata 105 (2004) 171--188.

\bibitem{Con00}
G.R. Conner,
Discreteness properties of translation numbers in solvable groups,
J. Group Theory 3 (2000) 77--94.

\bibitem{Deh98}
P.~Dehornoy,
Gaussian groups are torsion free,
J. Algebra 210 (1998) 291--297.

\bibitem{Deh02}
P.~Dehornoy,
Groupes de Garside,
Ann.{} Scient.{} Ec.{} Norm.{} Sup.{} 35 (2002) 267--306.

\bibitem{DP99}
P.~Dehornoy, L.~Paris,
Gaussian groups and Garside groups, two generalisations of Artin groups,
Proc. London Math. Soc. 79 (1999) 569--604.

\bibitem{EM94}
E.A.~Elrifai, H.R.~Morton,
Algorithms for positive braids,
Quart. J. Math. Oxford Ser. 45 (1994) 479--497.

\bibitem{Gar69}
F.A.~Garside,
The braid group and other groups,
Quart. J. Math. Oxford  Ser. 20 (1969) 235--254.

\bibitem{Geb05}
V.~Gebhardt,
A new approach to the conjugacy problem in Garside groups,
J. Algebra 292 (2005) 282--302.

\bibitem{GS91}
S.~Gersten, H.~Short,
Rational subgroups of biautomatic groups,
Ann.{} Math.{} 134 (1991) 125--158.

\bibitem{Gro87}
M.~Gromov,
Hyperbolic groups,
in: S.{} Gersten (Ed.),
Essays in group theory,
Math. Sci. Res. Inst. Publ., 8,
Springer, New York, 1987, pp. 75--263.


\bibitem{Kap97}
I.~Kapovich,
Small cancellation groups and translation numbers,
Trans.{} Amer.{} Math.{} Soc.{} 349 (1997) 1851--1875.

\bibitem{KR03}
D. Kim, D. Rolfsen,
An ordering for groups of pure braids and fibre-type hyperplane arrangements,
Canad. J. Math. 55 (2003) 822--838.

\bibitem{LL06a}
E.-K.~Lee, S.J.~Lee,
Some power of an element in a Garside group is conjugate to a periodically geodesic element,
arXiv:math.GN/0604144.

\bibitem{LL06b}
E.-K.~Lee, S.J.~Lee,
Abelian subgroups of Garside groups,
arXiv:math.GT/0609683, to appear in Comm. Algebra.

\bibitem{Lee07}
S.J.~Lee,
Garside groups are strongly translation discrete,
J.~Algebra 309 (2007) 594--609.


\bibitem{LM71}
S. Lipschutz, C. F. Miller III,
Groups with certain solvable and unsolvable decision problems,
Comm. Pure Appl. Math. 24 (1971) 7--15.


\bibitem{Pic01b}
M.~Picantin,
The conjugacy problem in small Gaussian groups,
Comm.{} Algebra 29 (2001) 1021--1039.


\bibitem{Sib02}
H.~Sibert,
Extraction of roots in Garside groups,
Comm.{} Algebra (2002) 2915--2927.

\bibitem{Sty78}
V.B.~Sty\v shnev,
Taking the root in the braid group (Russian),
Izv.{} Akad.{} Nauk SSSR Ser.{} Mat.{} 42 (1978) 1120--1131.

\bibitem{Swe95}
E.L.~Swenson,
Hyperbolic elements in negatively curved groups,
Geom.{} Dedicata 55 (1995) 199--210.

\bibitem{Thu92}
D.B.A.~Epstein, J.W.~Cannon, D.F.~Holt,
S.V.F.~Levy, M.S.~Paterson, W.P.~Thurston,
Word processing in groups,
Jones and Bartlett Publishers, Boston, MA, 1992, Chapter 9.
\end{thebibliography}
\end{document}